\newif\ifabstract
\abstracttrue
\abstractfalse 
\newif\iffull
\ifabstract \fullfalse \else \fulltrue \fi

\documentclass[11pt]{article}
\usepackage{amsfonts}
\usepackage{amssymb}
\usepackage{amstext}
\usepackage{amsmath}
\usepackage{xspace}
\usepackage{theorem}
\usepackage{graphicx}
\usepackage{url}
\usepackage{graphics}
\usepackage{colordvi}
\usepackage{colordvi}
\usepackage{subfigure}
\usepackage{lscape}
\usepackage{natbib}
\usepackage[english]{babel}

\advance \oddsidemargin by -0.8in
\newcommand{\myparskip}{3pt}
\parskip \myparskip

{\hspace*{\fill}$\Box$\par\vspace{4mm}}

\newtheorem{theorem}{Theorem}[section]
\newtheorem{corollary}{Corollary}[theorem]

\begin{document}

\title{Integrated Microsimulation Framework for Dynamic
Pedestrian Movement Estimation in Mobility Hub \footnote{An extended abstract appeared in IATBR 2015}}
\author{Alexis Pibrac \thanks{Laboratory of Innovations in Transportation (LITrans), Department of Civil, Geotechnical, and Mining Engineering, Polytechnique Montr\'eal, Montr\'eal, Canada, Email: {alexis.pibrac@polymtl.ca}} \and Bilal Farooq\thanks{Laboratory of Innovations in Transportation (LITrans), Department of Civil, Geotechnical, and Mining Engineering, Polytechnique Montr\'eal, Montr\'eal, Canada, Email: {bilal.farooq@polymtl.ca}}}

\begin{titlepage}
\maketitle

\thispagestyle{empty}

\begin{abstract}
We present an integrated microsimulation framework to estimate the pedestrian movement over time and space with limited data on directional counts. Using the activity-based approach, simulation can compute the overall demand and trajectory of each agent, which are in accordance with the available partial observations and are in response to the initial and evolving supply conditions and schedules. This simulation contains a chain of processes including: activities generation, decision point choices, and assignment. They are considered in an iteratively updating loop so that the simulation can dynamically correct its estimates of demand. A Markov chain is constructed for this loop. These considerations transform the problem into a convergence problem. A Metropolitan Hasting algorithm is then adapted to identify the optimal solution. This framework can be used to fill the lack of data or to model the reactions of demand to exogenous changes in the scenario. Finally, we present a case study on Montr\'eal Central Station, on which we tested the developed framework and calibrated the models. We then applied it to a possible future scenario for the same station.
\end{abstract}

\end{titlepage}

\section{Introduction}
\noindent
With the constant increase in the population of urban areas around the world, transportation and logistic is facing more organizational problems in order to deal with complex networks, mixing new technologies, and modern modes of transport. Never in the history, society has offered such a number of different possibilities, from the traditional individual modes (such as cars or bikes) to new concepts born in the growing market of sharing economy. Public transit systems such as metro, bus or tramway are now available in all sufficiently big cities. Moreover, these cities are well interconnected, thanks to various long distance modes of transportation. Thus rendering the current network of transportation facilities highly  efficient as well as highly complex. Despite the improvements in transportation technologies and increasing demand, the mode that has always remained central is the walking mode. Mobility hubs (e.g. train stations, terminals, etc.) within which walking is the only mode, are the key connections in the dominantly prevalent inter-modal urban travel patterns. They have a high risk of overcrowdedness and thus playing more and more prominent role in the fluidity and efficiency of the whole transportation network.\\

Despite significant advances in the individual level microscopic models to describe and reproduce pedestrians movement, the main limitation for simulations remains the lack of data for such situations. Indeed, with enough data, one particular case can be reproduced in a consistent manner (not exact but at least representative). But problems may arise when the information is incomplete; when it comes to validation (where additional data are needed for another time period); or extension of the scenario to future situations (where data are impossible to get).\\

That is why here, we develop a novel framework for pedestrian dynamics in which the demand part is no more a static estimation directly obtained from the data. The demand which is technically a time dependent Origin-Destination matrix, will for sure be based on the available data, but will also be influenced by other kinds of information, such as the schedule of transportation systems, infrastructure in which the agents are moving, estimates of the transfer times for each trajectory, etc. By bringing in new processes and dynamic supply information, we aim to account for incomplete data when it comes to generating the exact demand using a microsimulation. Furthermore, we will be able to estimate changes in this demand induced by the changes in exogenous inputs. For instance if the design of a train station has changed, we need to adapt the departure time of each individual following what would be their reaction in real life. In such a case, the demand still depends on the observations already gathered, but the link between them is no more direct. Some of these observations may not be exactly satisfied, but adapted depending on the changes we made in the scenario. Thereby now that the demand description is adapted depending on the situations, we are making a step forward in terms of realism, when it comes to simulating non-existent scenarios i.e. testing potential future changes. \\

Traditionally, the demand part of a scenario has been the starting point of a simulation--especially in case of pedestrian simulations \citep{Abdelghany2016, sahaleh2012scenario}. For example, in the four-step model, after the generation and distribution steps, the demand is completely described. Then comes the modal choice and assignment that are using the so-called demand and models like discrete choice theory, model of transport modes, etc. that describe the behavior of each agent. In this approach, the simulation is divided into two phases: first we create the demand, and then we use it into successive behavioral models. In fact, we create an agent and its characteristics, then we describe its movement thanks to a description of its behavior. And this behavior is simulated with a chain of models that are successively going deeper in term of information (first only the mode of transport is chosen, then the global itinerary is computed etc. until we obtain the complete time dependent description of the movement). Here the demand is no more considered completely exogenous or known a priori, but dependent on other parts of the scenario, that can also be partial results of behavioral models. We can't consider its generation into a separate phase. We have no longer a clear chain of objects to generate in a simple order thanks to deterministic models. But we have to find an equilibrium between all the different parts of the state, verifying all dependencies settled between them. The behaviors depend on the demand, for example the transfer times of each travel or the occupancy of each transport mode is directly influenced by the number of pedestrian in the station and their temporal distribution. And the demand depends on the behavioral simulation results, for example, the departure time is impacted by the time agents need to transfer or availability of modes.\\

\begin{figure}[!h]
	\centering
	\includegraphics[scale=0.43]{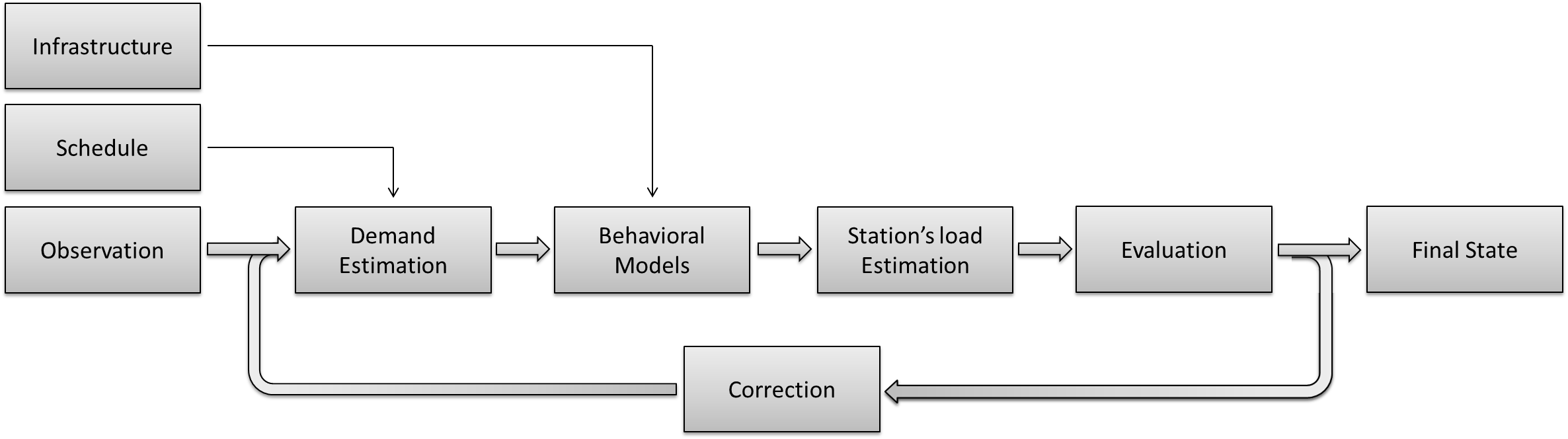}
	\caption{
		\label{markov}Organization of the proposed framework.}
\end{figure}

The problem of finding such an equilibrium is analogous to the one in Dynamic Traffic Assignment for vehicular traffic. However, due to the presence of a well-defined network and clear constraints, the search process for equilibrium in vehicular network is relatively trivial. Due to the complex movement of pedestrians and the high number of external factors that influence it (for example, arrival or departure time of a bus or a train is such an external factor that does not exist in vehicular traffic), the resolution for the case of pedestrian is of a higher complexity. To solve our equilibrium search problem, we are using a similar solution: a looping process running several times the same models until the convergence is reached. The classic behavior models of pedestrians simulations will be looped and computed as long as an equilibrium has not been found (see Figure \ref{markov}) i.e. until the simulated demand is consistent with the state generated from the partially observed demand. The purpose here is to present a novel microsimulation framework that controls the generation of the demand, intended movement patterns, and assignment in order to search for the equilibrium. In next section we present the existing work, after which the core methodology is presented. The case study of Montr\'eal Central Station is developed as an implementation of the proposed methodology. The results of base case and future scenario are discussed in details. In the end we present the conclusions and future direction.

\pagebreak
\section{Literature Review}
\noindent
Extensive research on various aspects of pedestrians simulation can be found in the literature. This has resulted in variety of tools to model the problem \citep{daamen2004modelling}. Past research has either focused on a specific operation within a train station \citep{zhang2008modeling} or the whole station \citep{sahaleh2012scenario}. The classical way to describe the pedestrian behavior is divided into three levels \citep{daamen2004modelling}: strategic, tactical and operational.
\\

The generation of OD matrix, that contains complete information about departure location, arrival location, and departure time for all agents, is a classical but tough problem in transportation research. It has been extensively studied in different contexts, e.g. for vehicles at urban area level for planning purposes \citep{national2012travel}, as well as at a smaller spatial scale like ours. Various available datasets have been used, beginning with traffic counts on the network in order to directly generate the matrix \citep{cascetta1988unified} or more recently with a Bayesian resolution \citep{cheng2014bayesian}. The schedule can also be used for this step  \citep{hanseler2015schedule}. These different information can be mixed in order to generate the matrix with a crucial time dependency \citep{ashok1996estimation}. Depending on the type of specific problem, a wide range of algorithms have been developed and tested, \cite{antoniou2014framework} provide an extensive literature review and propose a framework to compare them.
\\

The tactical level is the process that affects for each pedestrian their global route, depending on the OD matrix. In this step, we consider that all agents think in a graph-styled simplified network that represents the practical space and decision points. The classical formulation of this problem is the search of a Nash equilibrium \citep{wardrop1952road}. For vehicular simulation, the tactical level proceeds to the route choice of all agents \citep{bovy1990route}. Similar works have been developed for pedestrians \citep{hoogendoorn2004pedestrian}. However, in case of pedestrians we are of the view that it is behaviorally more consistent to consider this step as selection of decision points. The pedestrian choose their way through a succession of crucial decision points at a rather aggregate and abstract level. For example: which door, or coffee stand, etc. The process can be similar to the way finding algorithms for urban navigation that often use a graph representation of the network \citep{gaisbauer2008wayfinding}. Contrary to a continuous simulation of the trajectories where the space of possible solutions is also continuous, this level is characterized by discrete choices and so a finite number of possible configurations. The discrete choice theories have played a crucial role in the transportation research \citep{ben1985discrete} since they are related to different levels, such as modal choice \citep{hausman1978conditional}. Finally the proposed process strongly depends on a route cost function that should take into account the main phenomena, such as the travel time \citep{avineri2006impact} or even the perception of the facilities \citep{sisiopiku2003pedestrian}.
\\

The operational level goes one step further in terms of precision. Using the high level paths generated from previous process, it computes trajectories of all agents. Variety of models have been developed in this context. The more efficient are often aggregate model, where agents are gathered in order to consider the whole crowd like a flow \citep{hughes2002continuum}. This kind of approach can also be solved with a Cell Transmission Model, that discretize the space into cells \citep{daganzo1994cell}. \cite{hanseler2014macroscopic} have developed the cell transmission based model for pedestrians. The main advantages of these approaches are a quick simulation/enumeration time and a relatively good aggregate level precision for real crowd despite overly simplified assumptions. But in our case, we are interested in precise results with information on each pedestrian. As all the other levels are individual level, we want to maintain the consistency and disaggregation at operation level as well. We are interested in a microscopic scale. Several models have been developed at micro-scale, such as the use of discrete choices to model the next step of pedestrians\citep{robin2009specification} or an analogy with physical forces called the social force model \citep{helbing1995social}. In these models, it is always possible to go deeper in description to have better precision. Some studies have developed even more complicated description of agents taking into account for example the social or natural effects such as the use of field of view \citep{turner2002encoding}. These kind of agent-based model are now efficient on complex networks \citep{batty2003agent} and bring depth to the analysis.
\\

Once the estimation of all trajectories have been done, our goal is to authenticate the previous departure time and to correct them if needed. In the literature, algorithms have been developed that include a choice in the departure time generation \citep{de2002real}. The clear advantage is that it coincides more easily with real traffic conditions. Other algorithms try to deal with a real time correction of OD matrices \citep{bierlaire2004efficient}. But the new problem we are facing is that previous state estimated by the simulation step doesn't match any more with the new departure times. These simulations need to be recomputed. We now have a loop and need to find a convergence (Figure \ref{markov}). This problem is known as the Dynamic Traffic Assignment \citep{peeta2001foundations}. Some recent works proposed processes in order to solve this kind of problem \citep{nagel2012agent}.
\\

We propose a stochastic approach to solve this convergence problem. The output of each process will no longer be deterministic, but subject to probabilities as it has been proposed in \cite{daganzo1977stochastic}. The outputs we will now consider are probability distributions over the possible states space. In such case, a Bayesian resolution can be used \citep{maher1983inferences}. Specifically, we propose to consider the series of processes as a Markov Chain, using only the previous estimation of the state and giving back a new one, following a stochastic rule. The Monte Carlo algorithms therefore be used in order to identify the most probable states. One of such algorithm, Metropolis-Hasting algorithm \citep{hastings1970monte} has already be used for route choice set generation in a complex traffic network with high numbers of alternatives \citep{flotterod2011bayesian}.

\section{Methodology}
\subsection{Problem Statement}
\noindent
Given the infrastructure $I$, schedule of all modes of transportation $C$ and the location of different considered activities $A$, we are interested in estimating the state $S$ of the station that matches as good as possible to a set of incomplete observations $D$. A state $S$ contains the complete information of each pedestrian i.e. their activity chain $A^i$, their start and end location $(l_s^i,l_e^i)$, their starting time $t_{dep}^i$ and their exact trajectory $T^i:t\mapsto l$.\\
\[
(I,C,A,D) \mapsto S=(A^i,l_s^i,l_e^i,t_{dep}^i,T^i)_i
\]

Indeed, if all these information were contained in $D$, the proposed framework is obsolete. However, in reality $D$ is not sufficient to directly extract $S$: $D\neq \subset S$. Moreover, we may have to confront cases where $D$ was collected in a different scenario than the one in which $I$, $C$ and $A$ are defined. This happens when $(I,C,A)$ represents a non-existent scenario (for example possible perturbations of the reality or extensive changes that may occurs in the future). $D$ always corresponds to a scenario that has already happened i.e. base case. In such case, $D$ is still bringing a necessary amount of information, but they will not be directly considered as constraints for $S$. A necessary level of abstraction have to be brought to these observations: for example if a pedestrian $j$ is observed at a certain point of the time and space (this information $D_j$ is contained in $D$) it will not necessary be the case in $S$, the information could be transformed into $\overline{D_j}=$"$j$ \text{is taking bus} $b$". In $S$, $\overline{D_j}$ can bring pedestrian $j$ to have another trajectory if bus $b$ has a different departure time, $D_j$ is not satisfied. By calling $I_D$, $C_D$ and $A_D$ the respective infrastructure, schedule and activity of the scenario where $D$ was observed, we can write:
\[
(I=I_D,C=C_D,A=A_D) \Rightarrow D \subset S(I,C,A,D)
\]

\label{ovD}
By calling $\overline{D}$ the abstract information of $D$. This set is defined, when $D\not\subset S(I,C,A,D)$, such as it verifies:
\[
\overline{D} \subset S(I_D,C_D,A_D,D)
\]
\[
\overline{D} \subset S(I,C,A,D)
\]

In cases where $D\neq S$, it means that one or more estimated states may represent D. These estimated states will only differ on the part where we have the lack of information. Note that we are assuming that there always exists at least one state in the search space that can verify all our constraints. This assumption is reasonable in the case where the search space is well-defined and D has enough information. The goal of the simulation is to fill in the exact amount of information needed and thus choose one final state S*. We can't assure that there will be a unique state to which the convergence can bring us. Since it would mean that we have perfectly described all human phenomena that come into effect in the station. In fact we only want a representative of what could happen in reality, just a consistent case that allows us to understand the main phenomena in the station. We will be able to converge to several different and completely consistent solutions. But if we don't bring enough constraints, this space of possible final solutions will be oversized, we need to restrain it enough to have usable results. This is why constraints such as schedule dependence and behavioral models' consistency will be added in these cases where we don't have enough data.

\subsection{Inputs}
\noindent
Different kind of inputs will be considered. The four main ones are the infrastructure, schedule, activity list, and observed data: $(I, C,A, D)$.
\begin{itemize}
	\item Infrastructure $I$: Spatial description of the infrastructure. Mainly composed of a CAD design model of the studied station, available facilities, and the main entrances.
	\item Schedule $C$: List of all transportation modes, with their arrival or departure time, location and capacity.
	\item Considered activities $A$: Description of all activities available, inside as well as outside the station, for considered agent. It should contain the type, location and possibly the time at which it is available. In case of mobility hubs, the prime activity is to go from one mode to another, so in this paper we will model a unique activity for every pedestrian. However, the proposed methodology can easily be extended to include full activity chain modeling.
	\item Observations $D$: These data can be of various form. The less precise are aggregated counts on different point of the station, for example the number of people entering/exiting it per unit minute of the scenario. More precise data can be incorporated if they are available, for example observation of the exact time each pedestrian entered the station (or cross a specific point); information on the origin and destination of each travel; or even some local trajectories observed within the field of view of cameras in the station directly. A detailed discussion on the types of data commonly available on pedestrians in public spaces can be found in \cite{Farooq2016}.
\end{itemize}

The different behaviorial models used in the simulation are also inputs: different results can be obtained depending on the accuracy of each model and their consistency with reality. As for the observation $D$, the models can be considered as constraints. Indeed there are constraints on the kind of behavior agents can have. The final state $S$ will have to verify these constraints to be consistent. We can bring more constraints with more restrictive models: where possible agents movement are more precisely defined. 

\subsection{Simulation Processes}

Here we use the chosen behavioral models, and the inputs $(I,C,A)$ in order to simulate pedestrian agents moving in the station and obtain a description of state $S$. The three levels of simulation (strategic, tactical and operational) are respectively implemented with the activities generation, decision points choice, and assignment models.

\subsubsection{Activities Generation}
\label{gene}
The generation phase aims to estimate the demand. At the end of this process we obtain a part of $S$: the number of pedestrians, the activity chain of each one, their start and end locations, and their starting time:
\[
(A^i,l_s^i,l_e^i,t_{dep}^i)_i
\]
Since we only consider one type of activity (work), the model we choose here for the generation is Location Choice Model (LCM) that assign a destination for each pedestrian. This model, coupled with an estimation of the occupancy for every transportation mode and a description of the variability upon time, is sufficient to generate the demand. Estimating demand of a new scenario exactly corresponds to the calibration of location choice model. Thanks to the information contained in $D$, our framework corrects the demand until it is consistent with all parts of the scenario by calibrating this model. We can then use it for other scenarios, where at least one of $(I,C,A)$ is changed. Such calibrated model contains exactly $\overline{D}$ (see Section \ref{ovD}). The information in $D$ is absorbed in the form of a model to have the abstraction necessary to be generic for several different scenarios.

\subsubsection{Decision Points Choice}
The decision point level generates a global movement pattern for each pedestrian depending on their origin and destination. In this phase, the station is viewed as a simplified network representing all different paths. At each node or decision point, a pedestrian is confronted to a choice scenario. Pedestrian chooses one of the possible direction towards the destination. At this level there is no description of time, and the pedestrian is not considering other agents or obstacle. But several kind of information can be brought, from a simple estimation of the different transfer times on each link to real information of perception: signs, sized of corridor, light etc. In our case we are using a basic model for the sake of simplification in the simulation i.e. shortest path model, but random utility based choice models can be used.
\subsubsection{Assignment}
Finally the assignment uses all information generated in the previous phases to compute the exact time-dependent trajectory of each involved agent. It depends on their global routes defined by decision points; on their interactions with other agents and obstacles; and on different personal characteristics that may change their behavior in order to represent the diversity. Here we used the social force model \citep{helbing1995social} in the simulations.
\subsection{Convergence}
After the three previous processes have been executed, a state is obtained. But, like traditional simulations, the demand was generated before the state of the station was estimated. This demand could not have used some crucial information on the station's load such as transfer times or occupancy, because they were not observable yet. However, this demand may be in total adequacy with the obtained results. The first step is to observe this adequacy or not and measure what is not coherent. Then this information can be used in a correction process that will correct the previous estimation of the demand, now that more information are available. This correction process close a loop that we can be represented as a Markov process. We will then use the Simulated Annealing algorithm, a special case of Metropolis-Hasting sampler \citep{ross2013}, to make it converge to the desired state.
\subsubsection{Corrections}
\label{twoscenarios}
The correction process is required to correct the estimated demand. As we saw in Section \ref{gene}, a set of rules formulated in a model and applied in the generation step results in obtaining the estimated demand for a scenario. The correction process has two kind of possible actions that directly define the kind of simulation direction we are interested in:
\begin{itemize}
	\item \emph{Calibration of the demand.} First application of the simulation is to calibrate our model used in the generation step. This simulation is based on the available observation $D$ and generates the exact demand to satisfy it. It corresponds to the creation and calibration of $\overline{D}$ that absorbs all the information of $D$.
	\item \emph{Simulation of unknown scenario.} Second application supposes that the first one has already found the equilibrium point in order to calibrate the demand generation step. It means that $\overline{D}$ has been created. This second application uses it to simulate a new scenario without the availability of $D$.
\end{itemize}
The correction process for the first application is used to correct the location choice model: based on the difference observed between the current estimated state and $D$, it changes its rules $\overline{D}$ to try a different search point. Concretely, in our case, the probabilities of LCM are changed. After the observations are made, a correction is chosen depending on the lack or excess of people going into each kind of location. This correction can be focused on a specific location trying to increase or reduce the number of agent interested on it, or can be a mix of different changes. At each application of the correction process, since the choice is random any correction can be applied, but the probabilities are made set such that the correction has better chance to correct the observed difference in a right manner.\\

The correction process for the second application is simpler: the estimated state is analyzed and the abstract conditions of $\overline{D}$ (that are gathered in the generation step) are tested. If some are not satisfied, the behavior of corresponding agents are changed with a probability according to the correction they need. For example if some agents  miss their bus or train (that they should take according to the generation model) their starting time is corrected. The same principle can be applied if certain occupancy of a transportation mode needs to be reached by adding or deleting agents in the simulation.
\subsubsection{Markov Process}
\paragraph{Construction of the chain.}
The correction we just defined closes the loop. When applied on an estimated state, it gives us a new potential state with some probability of being chosen. This loop can be considered as the transition of a Markov chain (see figure \ref{markov}). In order to use the powerful properties of Markov chains, we have to prove that the one we defined is one. This is done here by proving the two properties: \emph{irreducibly} and \emph{aperiodicity}.\\

These properties are easily proved in the case of chains applied in finite space of state. This is definitely not the case here: the probabilities in LCM, that define the current position of the state in search space, are continuous (from 0 to 1). The space is not finite, neither discrete. In the case of such continuous space of state, it is common to use distribution probability (on which the chain is applied) in order to find the same properties than discrete spaces. In our case, it is not possible since the transition probabilities are not obtained with a formal function that we can easily integrate or apply in a region of states. Our transition is a function we can compute on only one state at a time. And since it is a whole simulation, we can't apply to a consequent number of state at each transition.

\begin{corollary}[Discrete Consideration]
	\label{consideration}
	In order to prove Markov chain properties, we are using another kind of consideration i.e. the location choice model can be continuous, but it is always applied on a finite number of pedestrians. There are infinite scenarios and demand that can be generated thanks to LCM, but for any particular scenario $(I,C,A)$, there is a maximum number of pedestrians that can be generated. In such case, the proportions of LCM may be continuous, but since they will be applied on a finite number, their effect is discrete. More precisely, around each LCM configuration of parameters, there is a small interval in which all other configuration of parameters have the same effect in a scenario. All these set of parameters correspond, in fact, to the same state. Finally we find that there is a finite number of different LCM in the particular scenario we are considering for the simulation.
\end{corollary}

\begin{theorem}[Irreducibility]
	We have to prove that there exists $N$ for which from any set of parameters $P$ we can reach any other one $P'$ in $N$ iterations with a non zero probability. A set has finite number of parameters in our LCM, so we can write $P=(p^j)_{j\in [1,J]}$, $P'=(p'^j)_{j\in [1,J]}$. In each iteration, at least one of the parameters is changed. The amplitude of this change has a maximum, let's call it $A$. $p^j\in [0,1]$ can reach any other parameter $p'^j$ in $|p^j-p'^j|/A < 1/A$ steps. The probability to jump from $p^j$ to $p'^j$ in $1/A$ steps is not zero since there is a finite number of values that can be reached--we just need to reach a value close enough to $p'^j$.
	
	There is $J$ parameters to change, each has a non-zero probability to be changed in any other value in $1/A$ steps. Moreover, each has a non zero probability to be chosen and changed at each iteration. It means from any state $P$, we can reach any other $P'$ in $J\times 1/A$ step with a non zero probability:
	\[
	N=\frac{J}{A}
	\]
\end{theorem}

\begin{theorem}[Aperiodicity]
	Corollary \ref{consideration} ensures that two "very close" sets of parameters can have the same effect in the generation process for one particular scenario. In fact it ensures that around one state $P$ there is a small open set (not empty) of parameter configurations that define the same state. Since a change at any iteration have a maximum amplitude of $A$ and minumum of $0$, the change made to a parameter can be so small that the new configuration of parameters is still in the open set of the same states. For any state, at any iteration, there is a non zero probability that we stay in the same state. The periodicity of all states can't be higher than 1. Our chain is aperiodic.
\end{theorem}

\subsubsection{Search Algorithm}
A Markov Chain Monte Carlo (MCMC) simulation process can be used to sample from the developed Markov chain. In particular Simulated Annealing algorithm is used to converge to an optimal state. Transition of the process is already defined, the algorithm need an objective function and a temperature to decide whether or not each new state will be kept.
\subsubsection{Objective Function}
The objective function drives the choice of search towards the optimal state. This state have to be consistent with all inputs we had: $(I,C,A,D)$ and the behavioral models. Even if we can integrate the behavioral models consistency in the objective function by implementing a rating system, we don't have to in our case. This can be easily done in future works. For example, it is possible to integrate the comfort (or security) appreciation for each pedestrian following a behavioral model that measure the perceived comfort of everyone. Integrating it in the objective function would lead to states where people tend to choose the travel by maximizing their evaluation of comfort.\\

\begin{figure}[!h]
	\centering \includegraphics[scale=1.2]{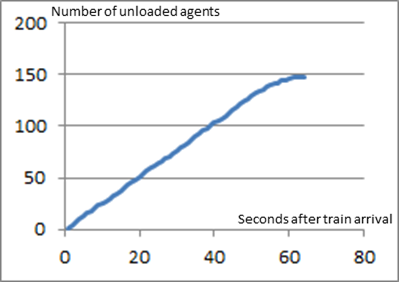}
	\qquad \includegraphics[scale=0.555]{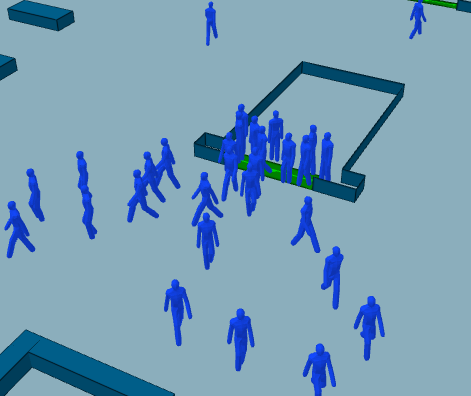}
	\caption{\label{UNLDM}Unloading model for trains: the bottleneck model is used due to the form of the connecting stairs between platforms and main hall.}
\end{figure}

The behavioral models are already used to generate the processes $I$ and $A$. Elements that the objective function should integrate to assure their impact on the simulation are $D$ and $C$: the conformity to external observation and the consistency with schedule, respectively.
\begin{itemize}
	\item Observations $D$: We use the comparison between these observation with the exact same information taken from the estimated state. $D$ is a set of values $(D_i)$ that correspond to a list of observation function $(O_i)$ applied on the real life station $S_{rl}$. These functions can be, for example, the number of pedestrians going through a particular door between two points in time.
	\[
	D=(D_i)=(O_i(S_{rl}))
	\]
	We evaluate these functions on the current estimated state $S$ to obtain $(s_i)$ the values to be compared with $D$. A rate is settled using the residual sum of square:
	\[
	OF_1(S,D)=RSS((s_i),(O_i(S_{rl})))=\sum_i (O_i(S)-D_i)^2
	\]
	\item Schedule $C$: the coherence of schedule measures the embarking and disembarking pattern for each train or any other mode of transport. It uses the list of pedestrians $(p_i^{(m)})$ taking each mode $m$. This list is obtained from the state $S$ thanks to the list of origin and destination of each pedestrian (time and space are considered) and the list of arrival and departure of each mode (time and platform also) in $C$. We assign a pedestrian going to or coming from a platform to a consistent bus or train.\\
	
	From this list of pedestrians we compute the arrival pattern $f(p_i^{(m)})$ of each mode for the state $S$. This pattern is compared to our embarking and disembarking pattern model $C_m$ (see Figure \ref{UNLDM}) and a rate is given with the residual sum of square to measure the consistency:
	\[
	OF_2(S,C)=\prod_m RSS(f(p_i^{(m)}),C_m)
	\]
	
	Finally, we may be unable to assign some pedestrians to a mode of transportation (origin or departure) if they are created before a mode arrives in the station or if they arrive too late to take the mode corresponding to their platform. We strongly penalize states with these incoherent observations. By denoting $Y(S,C)$ the number of such incoherent pedestrians in state $S$ with the schedule $C$, we have:
	\[
	OF_3(S,C)=e^{-Y(S,C)}
	\]
\end{itemize}
The relative importance of the three functions can be settled with two parameters $\alpha$ and $\beta$. We obtain the objective function:
\[
OF(S,C,D)=OF_1(S,D) \quad OF_2(S,C)^\alpha \quad OF_3(S,C)^\beta
\]

\pagebreak
\section{Implementation}

\begin{figure}[!h]
	\centering
	\includegraphics[scale=0.5]{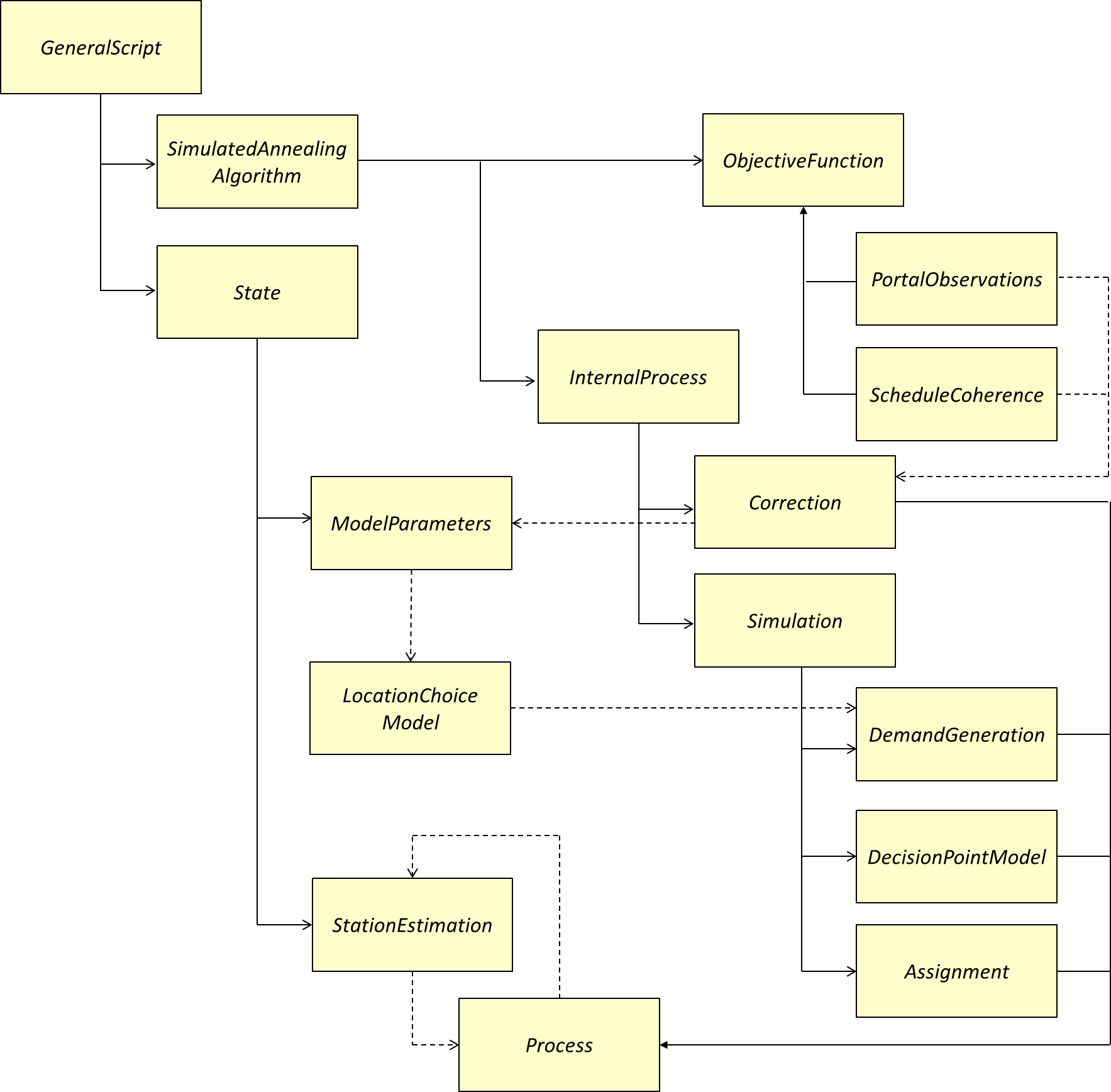}
	\caption{
		\label{imple}UML diagram of the framework.}
\end{figure}

Since the general algorithm, various behavioral models, and types of data are separate entities, we implement them in a way that each component can be independently and easily plugged in or replaced. We use an object oriented paradigm to implement in Java programming language. The implementation is available upon direct request to the corresponding author. Figure \ref{imple} shows the UML diagram of the framework. Different kind of scenarios and type of states can be plugged to the corresponding classes. The implemented code can be used in many different scenarios and is able to take various kind of models. Please also note that a commercial software called MassMotion by Oasis Software is used for running the \emph{Assignment} process.

\section{Case Study}
\noindent
As a case study we explored our framework on the Montr\'eal Central Train Station. Here 14 tracks are exploited by several national and local railways companies. The station is also linked to two metro stations, 16 bus lines, encloses an active underground mall, and is directly connected to several buildings. It is an important part of the Montr\'eal city centre since it is located downtown, and is a central part of the Montreal's Underground City, the biggest pedestrian indoor network in the world.
\subsection{Presentation}
\subsubsection{Simulation Setup}
\begin{figure}[!h]
	\centering
	\includegraphics[scale=0.7]{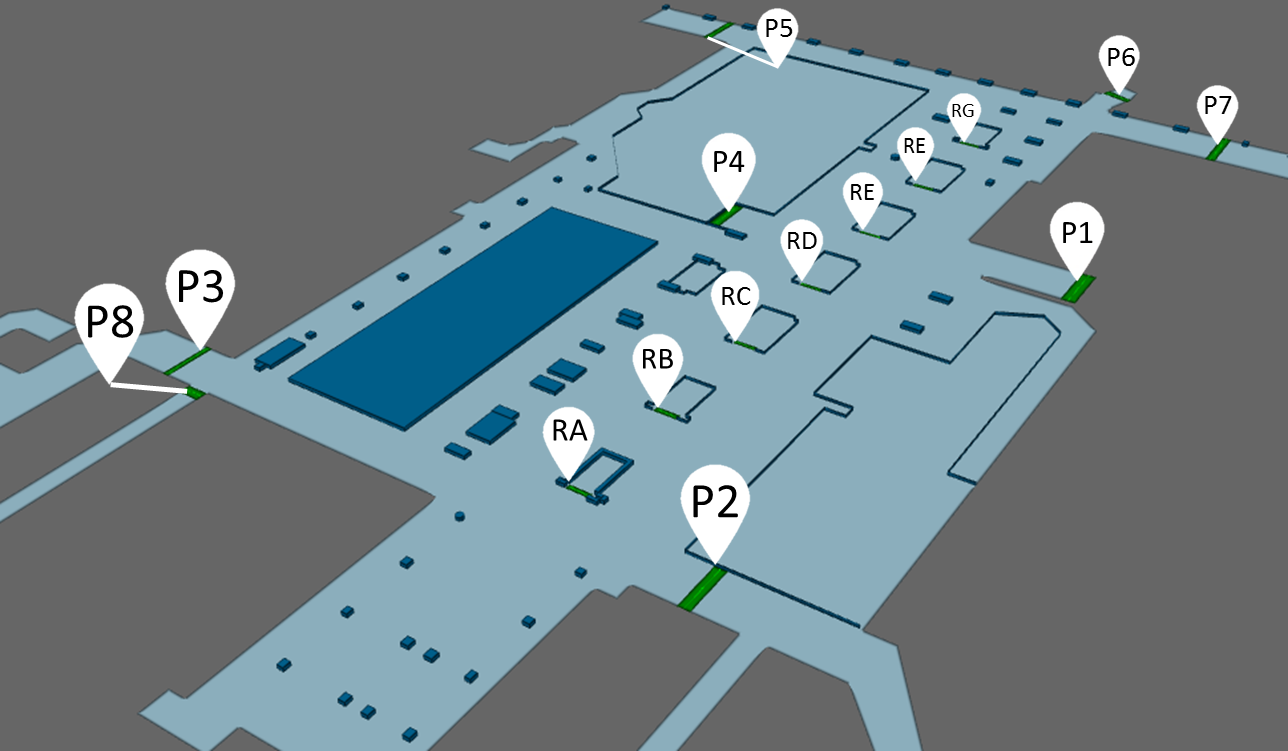}
	\caption{\label{station} Representation of the station. Usable space is in light blue and obstacles in dark blue. Portal with which people can enter and leave the simulation are in green.}
\end{figure}

\begin{figure}[!h]
	\centering
	\includegraphics[scale=0.7]{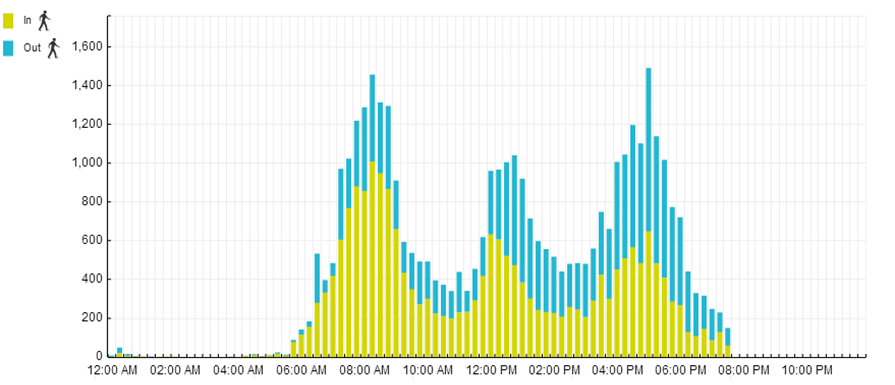}
	\caption{\label{data}Data recorded on October 1, 2014, at portal 2. Pedestrians using the portal as entrance are in green and leaving through it are in blue.}
\end{figure}

\noindent
We will study a fixed part of the space and time of the scenario i.e. we will model the main hall of the central station (see Figure \ref{station}). Pedestrians will be able to enter and leave through different portals that model the entrance/exit of boundaries of the station. These portals are the different corridors (1 to 8) arriving into the main hall of the station, and also the stairs connecting the platforms just under the hall (RA to RG), where trains arrive. The time of day that we are interested in is when the station is the most crowded i.e. the peak period. Since the afternoon peak period is more spread, it is less intense (see Figure \ref{data}). We have chosen the morning peak period. Agent based simulations are computationally very demanding and because we are running several simulations in a single iteration, we began with a short window of time (i.e. 15 minutes of highest demand) to minimize the computational time. So our simulation concerns the duration between 8:30am and 8:45am, the most crowded quarter, during which several trains arrive and leave the station.

\noindent
The schedule $C$ gathers all departure and arrival trains of a normal day, including their capacity information. During the time window we are studying, several of them are unloading and others are taking passengers for the suburb.

\subsubsection{Behavioral Models}
\label{models}
\noindent
For this first simulation we selected basic behavioral models. The three levels have to be implemented with one model: strategical, tactical and operational. In the strategical level, we should model the activity chain. Because we are simulating the morning peak hour, we assume that the main purpose of the displacements is work. The only dimension to generate here is the location of this activity. That is why we use an activity location choice model. This is particularly consistent because the studied space (the hall of central station) is small and doesn't host too many different activities. There are still some coffees and restaurants. In future simulations the model could integrate them and propose a full activity chain modeling.\\

The operational level, where we model decision points, is also impacted by the size of the station i.e. when there are not too many different ways to go from one point to another, its importance is diminished. We choose the simplest model there is, the shortest path. It is still particularly consistent since everyone is going to work at that time, people may mainly choose their trajectory to go as fast as possible. Finally the operational level is very important in term of realism. We used the social force model that gives a good description of real behavior. The parameter values calibrated in \cite{sahaleh2012scenario} were used.

\subsubsection{Scenario Development}
As explained in Section \ref{twoscenarios}, two kinds of simulation are possible: one using observations $D$ on existing use so as to calibrate the demand generator model. The second kind of simulation use the calibrated model to estimate the station's load in a scenario for which $D$ does not exist. For this case study we are executing both kind of simulations. First, thanks to the data we collected on the real scenario, we will calibrate the behavioral models. Then it will be used in a possible future scenario for which no data could be collected.

\subsection{Base Case Scenario}
$I$, $C$ and $A$ are already detailed, as well as the behavioral models. Only $D$ is needed to launch the simulation.
\subsubsection{Inputs}
\noindent
We gathered observations by installing magnetic sensors on all entrances of the station's main hall i.e. portal 1-8. Unfortunately, logistical issues bared us to measure the flow at access points to platforms. We also did not have access to the occupancy data of trains. The commercial sensors were provided by the manufacturer, Eco-Counter. Due to the data collection rate of these sensors, pedestrian counts were only recoded for 15 minute intervals during one typical day. Note that the pedestrian loading at portals in the simulation was done at 1 minute interval. So for this purpose the initial departures at entrance portals 1-8 were assigned based on Poisson Arrival Process for every 15 minutes of sensor counts. The departures from portals connected to platform were based on the arrival times of the trains and unloading curve.

\subsubsection{Optimal Solution Search}

\begin{figure}[!h]
	\centering
	\includegraphics[scale=0.63]{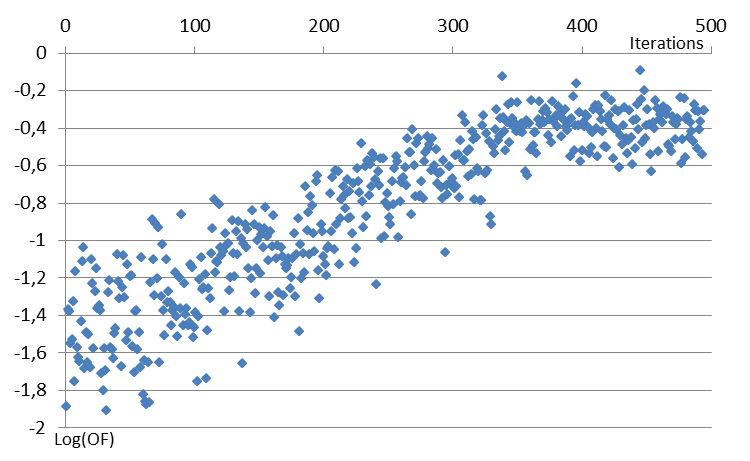}
	\includegraphics[scale=0.63]{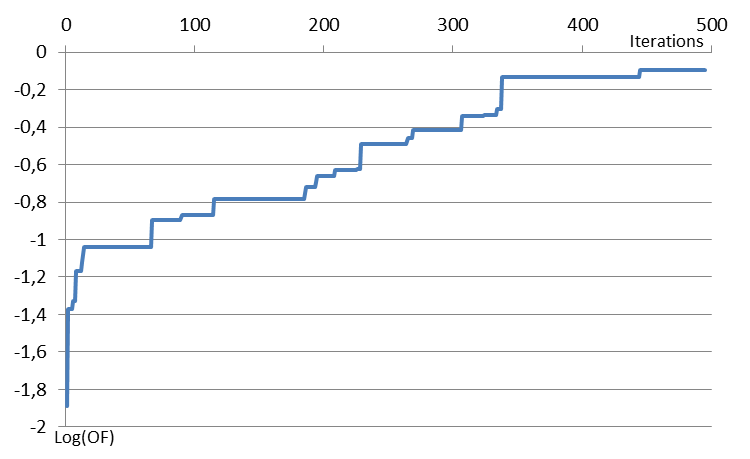}
	\caption{\label{c500}Left: $log(OF(S_i))$ for the state generated at each iteration $i$. Right: $log(OF(S^*))$, value of the objective function at step $i$ i.e. its value for the best known state.}
\end{figure}

\begin{figure}[!h]
	\centering
	\includegraphics[scale=0.5]{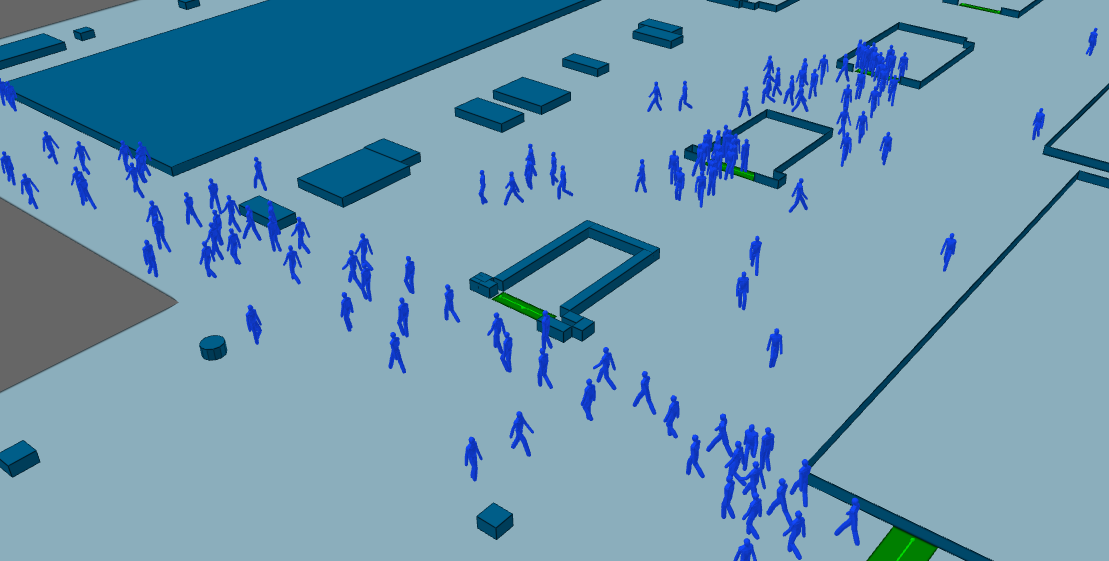}
	\caption{\label{screen}3D representation of a simulation. Two trains just arrived. The disembarking passengers mix with a continuous flow of pedestrians crossing the station.}
\end{figure}

After 500 iterations, the simulation converged to an optimal solution. Figure \ref{c500} shows values of the objective function for each iteration, and for the selected states. We can observe that there is a gradual and steady progression towards search regions with better values and thus the selected value constantly improved towards optimal solution. The convergence is very slow{---}it took several hundreds of iteration to obtain an acceptable result. The reason is that, for this particular simulation, we begin with a particularly incoherent set of parameters for the location choice model. The goal here is to show the robustness of the method i.e. that it converges, though slowly, to a proper solution. For sure, when the goal of a simulation is only to have consistent results, we can begin with a more coherent set of parameters, simply generated with the common sense of the analyst. Final solution is an estimation of the station's load with the trajectory of all pedestrians. For illustration purposes, we can see a 3D representation of the pedestrian movement in Figure \ref{screen}.

\subsubsection{Validation}

In order to validate the results, we need real observations that have not been used in the convergence process. The problem is that the lack of data is exactly what we try to solve here. So, instead, we used real train occupancy information in order to validate the coherence with the real life. Since people disembarking trains can take several other exits than through the main hall, we particularly compare the count of pedestrian leaving through the ones we considered. In the simulation, we found that an average of 760 people using these exits after the arrival of a train. In the real data we have an average of 850 people disembarking from the trains, which is higher. This difference can be explained because some exits from platform to the hall were not considered, so the flow is limited, in the simulation, to the principal exits only.\\

\begin{figure}[!h]
	\centering
	\includegraphics[scale=0.85]{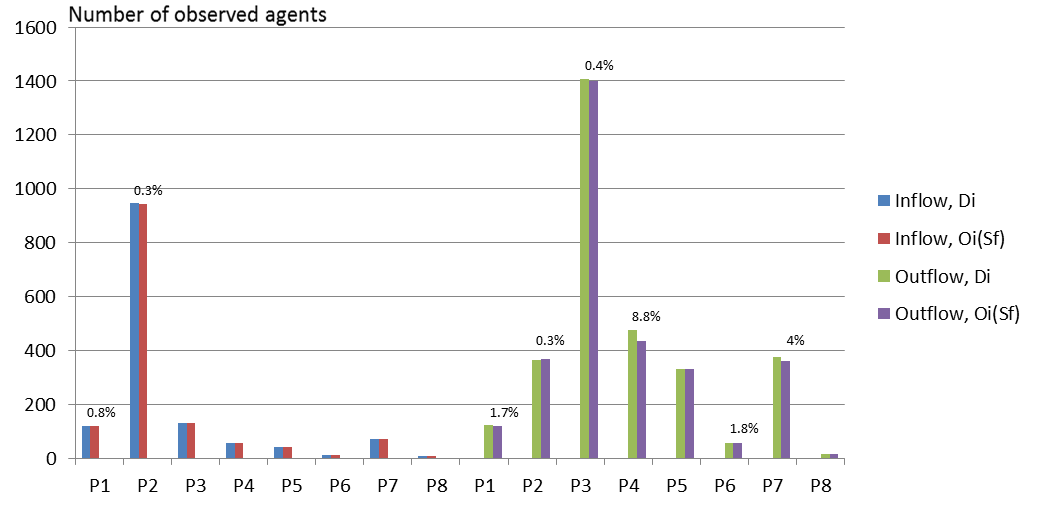}
	\caption{\label{table}Comparison between the real life observations ($D_i$) and the same observation on the final state ($O_i(S_f)$) for the inflow and outflow of all portals of the station. The percentage of error is written when it is not 0\%.}
\end{figure}

In order to validate the convergence, we can also analyze at what point the observed demand is satisfied by the solution. Figure \ref{table} shows fit of the optimal solution with the observed data. We can clearly see that most of conditions are satisfied with only one exception. According to the observation, more people should be leaving the hall through portal 4, but this error is less than 10 \%.\\

\subsection{Future Scenario}
After calibrating the simulation, we used it in a future scenario for which we did not have any real observations, but the scenario close enough to base case. Thus the utilization of base case calibrated models was consistent. We simulated here the same station, with the exact same facilities and people using it, but with an increase in the population by 50\%. This is a possible and very realistic scenario, if the infrastructure at station is not updated in the near future.
\subsubsection{Inputs}
The inputs are almost the same: $I$, $A$ and $C$ are unchanged, as well as the behavior models. Only $D$ is no more used. The LCM is now simply used instead of being calibrated. The total number of agents involved in the station is multiplied by 150\%.
\subsubsection{Simulation}
Even if the LCM is now static, we still need to make the Markov Chain converge. The demand still has to find an equilibrium with the estimation of the station. For example, an estimation of the transfer time is used in the generation step so that, after using the LCM, a departure time is assigned for each pedestrian. This is after several iterations that these estimations are consistent with the scenario so that the demand is generated in a consistent way.
\subsubsection{Results}
\begin{figure}[!h]
	\centering
	\includegraphics[scale=0.5]{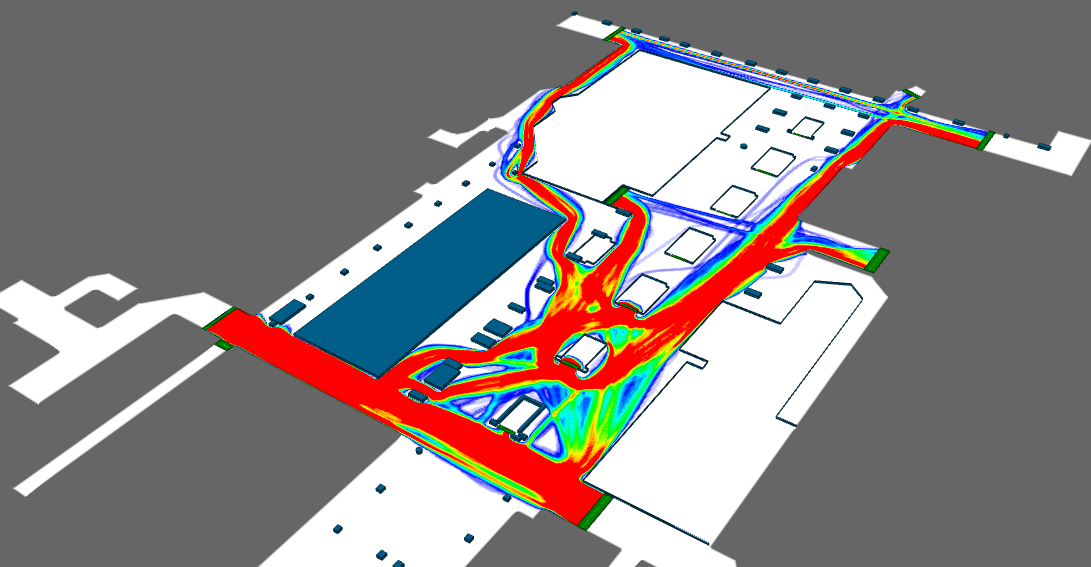}
	\caption{\label{flowR} Densities of pedestrians on each path}
\end{figure}

\begin{figure}[!h]
	\centering
	\includegraphics[scale=0.6]{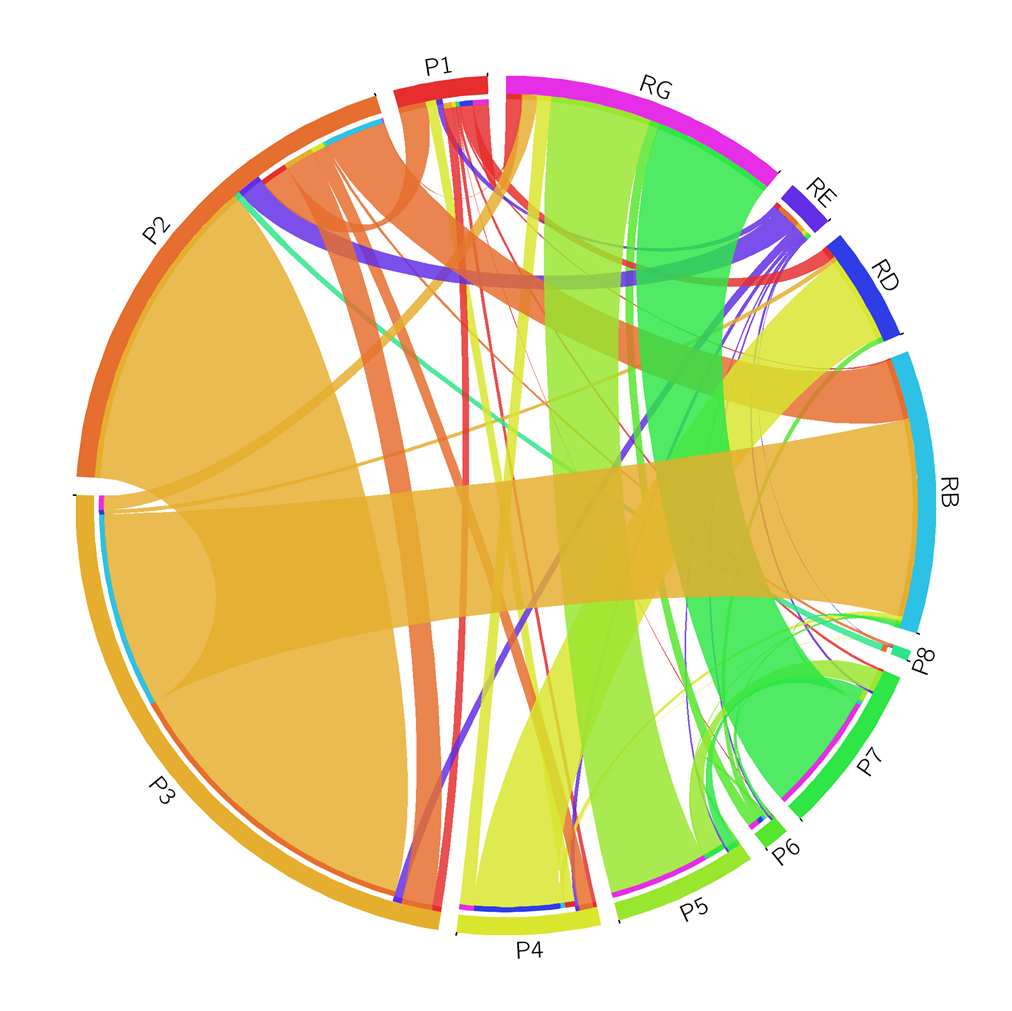}
	\caption{\label{flowL} Flow between each portal of the main hall.}
\end{figure}

From the converged state, we can extract information on pedestrian trajectories over space and time. Figure \ref{flowR} shows the principal paths used by pedestrians during the simulated time window. We can see what parts of the station are overcrowded and may present a risk of traffic congestion. The results also provide information on each agent. For example, a criterion could be used to measure the safety or satisfaction of each agent. The general OD matrix over time in the station can also be obtained from the pedestrian trajectories in the simulation. We represented it in Figure  \ref{flowL}. Each strip represent a flow from an origin (to which it is attached) to a destination. Its thickness is proportional to the number of agent using it. We can identify that the major flow is from portal 2 to 3. The two most loaded trains are arriving in platform B and G. Passengers arriving with the first one are mainly going to portal 2 and 3, while those arriving with second are oriented to portal 5 and 7. Only one platform is considered as a destination by pedestrians i.e. platform E. This is consistent since it is from this platform that the only train leaving in our simulation window time departs.\\

The comparison between base case and possible future scenario can bring a detailed picture of the evolution of station. Figure \ref{density} shows the densities over time in both scenarios. Densities are represented according to the standard \cite{fluin1971} and IATA (International Air Transport Association) level of service mappings. We can clearly see how an augmentation of the population in the station does not linearly increase the measured densities. With the 50\% augmentation, the presence of higher densities explodes. Indeed, when serious congestion appears, pedestrians get blocked and stay longer in the station. This leads to even more pedestrians in the station and so a higher danger of congestion--it a vicious circle phenomenon. Also we can identify some details of what parts could require some improvements. We can see in Figure \ref{congested} the different intersections where high levels of congestion appear.\\

\begin{landscape}
	\begin{figure}[!h]
		\centering
		\includegraphics[scale=0.3]{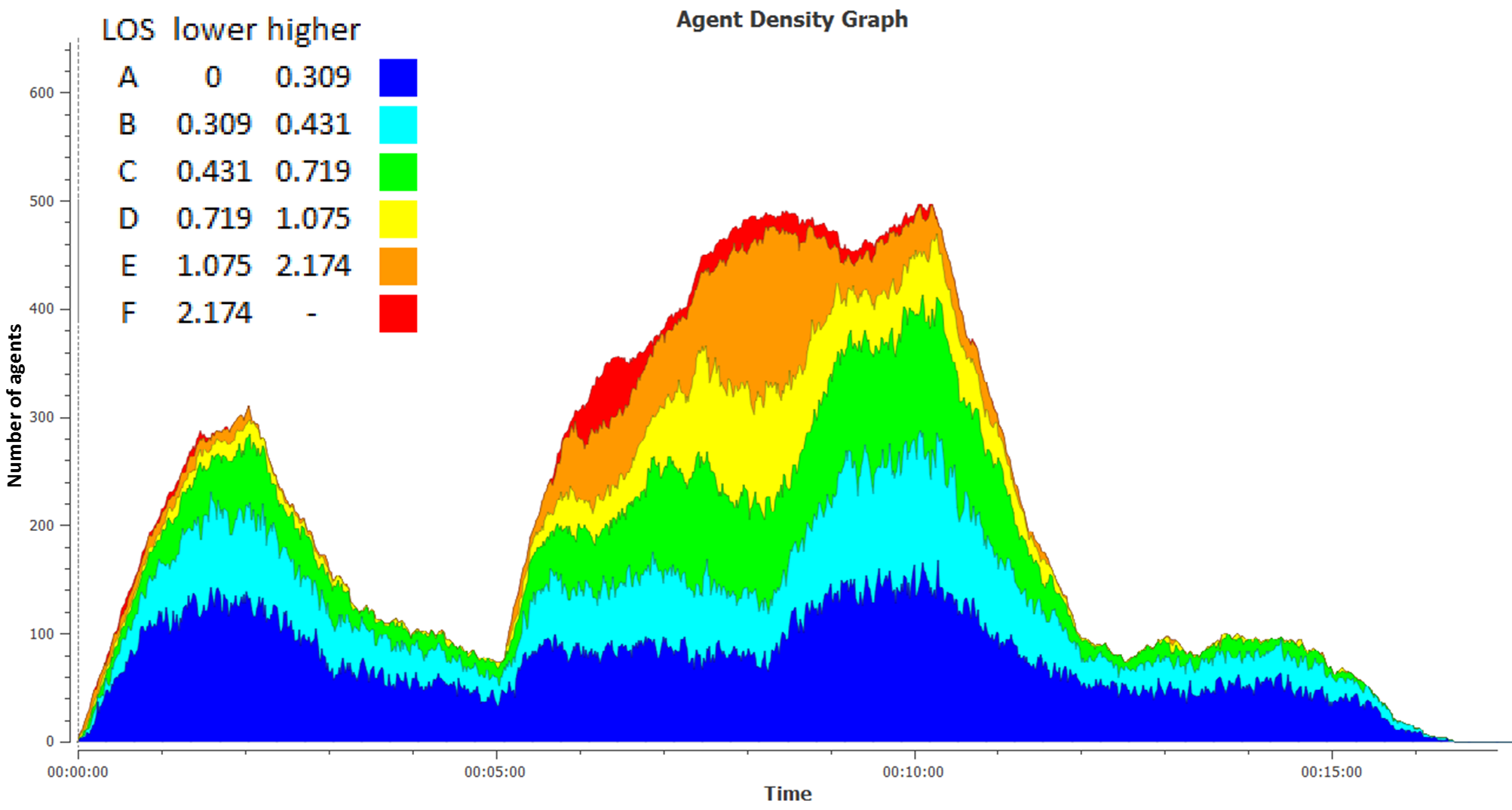}
		\includegraphics[scale=0.3]{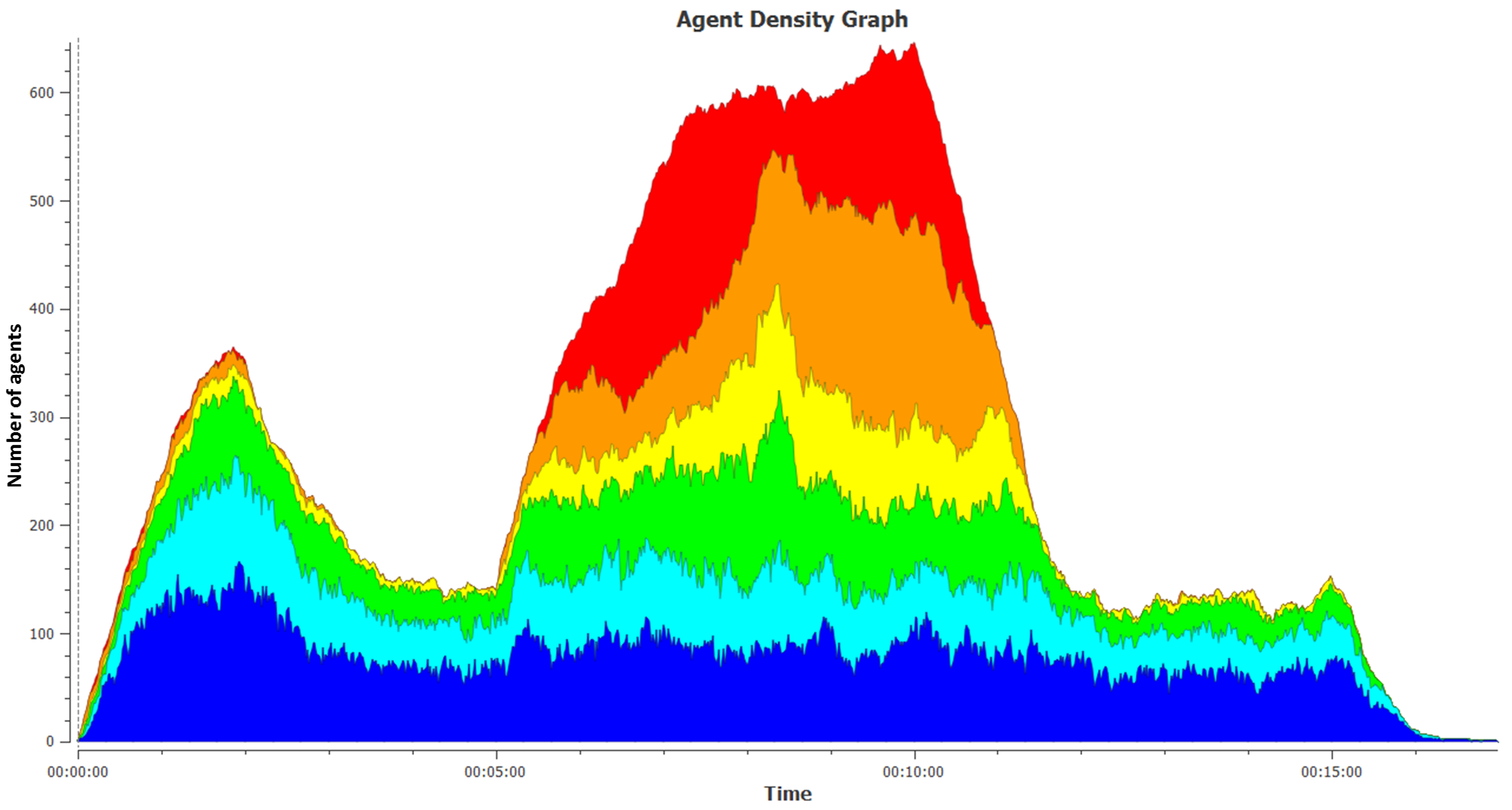}
		\caption{\label{density}Representation of densities over time in the station for base case (left) and future scenario (right).}
	\end{figure}
	
	\begin{figure}[!h]
		\centering
		\includegraphics[scale=0.87]{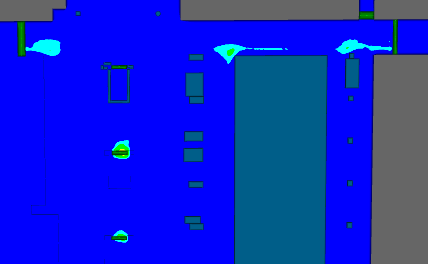}
		\includegraphics[scale=0.92]{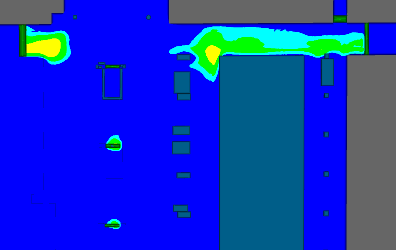}
		\caption{\label{congested}Spatial representation of densities for base case (right) and future scenario (left).}
	\end{figure}
\end{landscape}

\subsection{Discussion}

The first point that we would like to discuss is the 8.8\% error with respect to outflow observations at portal 4 (see Figure \ref{table}). This error means that our model did not manage to send enough people to this portal. The Location Choice Model is responsible for this assignment. The error corresponds to a default aspect of our model, which defined a type of attraction for each portal: city, metro, train, etc. The probability for pedestrians to choose one of these attractions was calibrated in the base scenario. Once an attraction was set for an agent, a destination was assigned by selecting the nearest portal that proposes this attraction. The lack of people going to portal 4 means that the attraction we assigned to this portal put it in competition with other portals that were surely closer to the major inflows. We can see in Figure \ref{flowL} that pedestrians leaving through portal 4 were essentially coming from platform D. And the flow created by this platform is limited.\\

The error means that the function of portal 4 was not properly assessed. In order to improve the model, we can use the random utility theory and define a utility function for the attractivity of each portal based on their attributes. We could also imagine to simply extend the location choice model by adding a type of attraction, just for portal 4. But we have to be careful with these options, since they mean more parameters to calibrate in the model. The search algorithm will be faced with a higher degree of complexity in search space to explore. More information will be needed for the algorithm to be able to determine an optimal solution.\\

We can identify in Figure \ref{flowR} another limitation of the implementation: some path used are not coherent with the reality. They cross an area of shops that is not attractive for pedestrians, in reality they try to avoid it and mostly take the wider corridor just next to the area. This difference between the real observation and the simulation is also coming from the model at tactical level i.e. the shortest path model. We observed that in the simulation, pedestrians are taking the path through shops because it is shorter and so it corresponds to the model. But in real life, the choice of path is more complex than a simple shortest path choice. When it comes to choosing between the two directions, people tend to take the corridor because it seems more attractive. Such phenomena in the choice are not described by the model. In future, we can suggest the use of a random utility based decision points model--especially the dynamic mixed logit model, which fits very well with the choice scenario. Such model could describe some human behavior such as the impact of the perceived environment in the choice, same person making successive decisions, and correlation between the decisions. People may tend to be more attracted by bigger corridors, shown by direction signs, and presenting less obstacles in the sides (such as tables or shops' advertising).\\



\section{Conclusion}

We presented an agent based microsimulation framework for pedestrian movement in moblity hubs and public spaces. The problem is formulated as a Markov chain of activity generation, decision points choice, and assignment processes. Thanks to behavioral considerations of the demand and dependence on public transit schedules, the resulting framework is truly dynamic and can fill the lack of complete observations. We propose MCMC process that converges to an optimal solution depending on the type of behvaioral models, infrastructure data, public transit schedules, and incomplete observed demand. As a result the framework is able to predict the activities, location, start time, duration, and detailed trajectory of individual pedestrian.\\

A case study of Montr\'eal Central Station has been implemented for the base case and a future scenario with demand augmentation of 50\%. The validation of the base case shows a good fit. We also observed several differences between the result and the data in station. They were all explained by the choices of model: a too simplistic description of the infrastructure and of the possible activities in the station; a calibrated location choice model not perfectly adaptable; and a shortest path model at tactical level that needs to be more representative of behavior and dynamic conditions. These are the dimensions where improvements can be made to the current implementation of the case study. The general algorithm is computed in a way that these changes can be easily integrated.\\

Finally, there is a great potential and applicability of the proposed microsimulation framework. Once behaviorally richer models are implemented, the simulation will be able to render how the station will be affected if specific changes are made to the design or schedule, with a dynamic demand that effectively reacts to these changes. For example a change in train's departure time will force pedestrians to leave at a different time in order to have a coherent behavior in the simulation. If the demand were not dynamic, these pedestrians' departure time could not be changed and we could observe absurd situations where the arrival time of the pedestrian is completely not coherent with his/her train. Such a framework will be very useful in the network-level optimization of the schedule for various modes of transportation, in order to have perfect connections between them following what the population needs, and avoiding high densities that could lead to unstable situations.

\section*{Acknowledgements}
We would like to thank Natural Sciences and Engineering Research Council of Canada (NSERC) and Fonds de recherche du Qu\'ebec Nature et technologies (FRQNT) for funding this research. We would also like to thank Soci\'et\'e de transport de Montr\'eal (STM), Agence m\'etropolitaine de transport (AMT), Oasis Software, and Eco-Counter for providing us the critical support to make the research possible in this paper.

\bibliographystyle{plainnat}
\bibliography{PibracFarooq_PedDynamics}

\end{document}